\documentclass[12pt]{article}
\usepackage{amsmath,amssymb,amsbsy,epsfig,amsfonts,graphics}
\newtheorem{Theorem}{Theorem}

\newtheorem{prop}[Theorem]{Proposition}
\newtheorem{cor}[Theorem]{Corollary}
\newtheorem{Lemma}[Theorem]{Lemma}

\newcommand{\proof}{{\bf Proof: }}
\newcommand{\Def}{{

\noindent\bf Definition: }}
\newtheorem{Corollary}[Theorem]{Corollary}

\newcommand{\qed}{\hfill\framebox(6,6){}}
\newcommand{\comment}[1]{} 

\newcommand{\A}{{\cal A}}

\newcommand{\F}{{\cal F}}

\newcommand{\Sig}{\Sigma}
\newcommand{\Del}{\Delta}
\newcommand{\del}{\delta}
\newcommand{\Lam}{\Lambda}
\newcommand{\lam}{\lambda}
\newcommand{\T}{{\cal T}}
\newcommand{\bound}{\partial}

\newcommand{\bbb}{\mathbb}

\def\<{\langle} 
\def\>{\rangle}

\newcommand{\Sum}{\sum}

\newcommand{\cF}{{\cal F}}
\newcommand{\diag}{\mbox{Diag}}
\newcommand{\ignore}[1]{}
\title{Affine $\Lambda$ Buildings II:\\
A reduction of axioms}
\author{Curtis D. Bennett\\ Loyola Marymount University}

\begin{document}

\maketitle

\section{Introduction}

Jacques Tits introduced the notion of a building as a geometry associated to groups of Lie type in \cite{T2}, providing new geometries associated to the exceptional groups of Lie type.   In 1972, F.~ Bruhat and Tits \cite{BT} developed a theory of affine buildings for the purpose of studying groups over fields having a discrete valuation, although their work applied more generally to groups over fields having a valuation over the real numbers.   Affine $\Lambda$-buildings were first introduced by the author in \cite{B1} and \cite{B2} as generalizations to the Bruhat-Tits buildings, allowing for groups over fields having a valuation into any totally ordered abelian group $\Lambda$ (and also generalizing the notion of a $\Lambda$-tree).   Recently, Linus Kramer and Katrin Tent have made use of affine $\Lambda$-buildings in their study of asymptotic cones and their short proof of the Margulis conjecture \cite{KT}, \cite{KST}.  

In \cite{B2}, the author defines an affine $\Lambda$ building as a pair $(\Delta,\cF)$ satisfying a set of six axioms given in section~2.  The first four of theses axioms are relatively easy  to check in most cases.  The sixth axiom is a little less straightforward and the sixth axiom can be particularly difficult to show.  Consequently, to make affine $\Lambda$-buildings more useful as a tool, it is worthwhile to find an easier set of axioms with which to work.   The purpose of this paper is to provide an easier set of axioms by extending results of Anne Parreau \cite{P}  on the equivalence of axioms for Euclidean buildings.  In particular, in section~3 we define a strong exchange condition, mimicking Proposition 3.27 from \cite{T3}.  In Section~4 we then prove that a pair $(\Delta,\cF)$ is an affine $\Lambda$-building if and only if it satisfies the first four axioms together with the strong exchange condition.   More recently, Petra Schwer has used these results to obtain  an even stronger conclusion \cite{H}.

\section{Preliminaries}

In the first part of this section we follow \cite{B2} to give the definition of an affine 
$\Lambda$-building.  In the second part of this section we give precise definitions of our 
new replacement axioms. 

\subsection{The $\Lambda$-apartment $\Sigma$}

Let $\Lambda$ denote a totally ordered Abelian group viewed as a $\bbb{Q}$-module.  Let $\overline{W}$ be a spherical Coxeter group, $\Phi$ the associated spherical root 
system with base $\{\alpha_1,\dots,\alpha_n\}$, and let $A=(\<\alpha_i,\alpha_j\>_\Phi)$ 
be the associated Cartan matrix, and let $D=\diag(d_1,\dots,d_n)$, be given so that $D^{-1}A$ is symmetric.  Define 
$$ 
\Sigma = \left\{ \Sum_{i=1}^n \lambda_i\alpha_i\mid \lambda_i\in\Lambda\right\} \cong \Lambda^n.
$$
We represent elements of $\Sigma$ by the coordinates $(\lam_1,\dots,\lam_n)$, and note that any sum of roots with coefficients in $\Lambda$ corresponds to a unique element of $\Sigma$.  Moreover, using this, we extend the action of $\overline{W}$ to $\Sigma$ linearly.  Letting $\T$ be the set of $\Lam$-translations of $\Sigma$ normalized by $\overline{W}$, we define the {\it affine $\Lambda$-Weyl group\/} to be the group $W=\T\overline{W}$.  Reflections of $W$ are defined to be conjugations of reflections of $\overline{W}$.

Walls of $\Sigma$ are defined to be the fixed points of reflections of $W$.  For $v\in\Sigma$ and a wall $M$ associated to the reflection over $\alpha_i$, using Proposition~2.1 of \cite{B2}, we have $v$ can be written uniquely as $v=m+\lam^i\alpha_i$ for some $m\in M$ and $\lam\in\Lambda$.  For each $i=1,\dots,n$, we let $v^i=d_i\lam^i$, and in general, we let $v^{w(\alpha_i)}=(w^{-1}(v))^i$ for $w\in\overline{W}$ (this is well defined by Corollary~2.5 of \cite{B2}).  As is shown in \cite{B2}, $x\in\Sig$ is uniquely determined by the $n$-tuple $(x^i,\dots,x^n)$.

A {\em half-apartment} of $\Sig$ is any of the sets
$$
H_w=w(\{v\in\Sig\mid v^i\ge0\}),
$$
and the {\em fundamental sector} of $\Sig$ is the set
$$
\widehat{S}=\{v\in\Sig\mid v^i\ge0\quad \mbox{for all $i=1,\dots,n$}\}.
$$
A set $S$ is a {\em sector} of $\Sig$ if $S=w(\widehat{S})$ for some $w\in W$, and we say that $S$ is based at $w(o)$ where $o=(0,\dots,0)\in \Sig$.  The {\em fundamental sector panel of type $i$} of $\widehat{S}$ is the set
$$
P_i=\{v\in\widehat{S}\mid v^i=0\},
$$
and we again use $W$ to define the sector panels of type $i$ of any sector.  

We now define a $\Lam$-distance function on $\Sigma$ as follows.  Let $\alpha_i$ be a fundamental root of $\Phi$ and $v=(\lam_1,\dots,\lam_n)$.  Define
$$
(\alpha_i,v)=d_ia_{ij}\lam_j
$$
where $A=(a_{ij})$ is the Cartan matrix associated to $\Phi$, and $D=\diag(d_1,\dots,d_n)$ is the given ``symmetrization'' factor of $A$ (as above).  For $\alpha\in\Phi$, we define $(\alpha,v)$ by extending the function linearly.  The $\Lambda$-distance between vectors $v_1,v_2\in\Sig$ is given by
$$
d(v_1,v_2)=\sum_{\alpha\in\Phi_+}|(\alpha,v_1-v_2)|,
$$
where $|\ |$ denotes the usual absolute value function on $\Lam$.  As shown in \cite{B2}, $d$ is a symmetric $W$-invariant function that satisfies the triangle inequality, and so it is appropriate to call $d$ a  $\Lambda$-metric.

With these definitions, we define two sectors to be {\em parallel} if they are at bounded distance from each other.  (Note that it follows that sectors are parallel if and only if their intersection contains a subsector.)  We similarly define parallelism of sector-panels, and in \cite{B2} it is shown that a sector or sector-panel $X_1$ is parallel to $X_2$ if and only if $X_1$ and $X_2$ are translations of each other. 

A set $\Omega\subset\Sig$ is {\em convex} if it is the intersection of half-apartments, and it is {\em closed and convex} if it is the intersection of finitely many half-apartments.  

\ignore{Define sector-germ, sector-gallery, and sector distance.}

\Def  Given a sector $S$ of $\Sig$ based at $y$, we say that a subset $V$ of $\Sigma$ 
contains the {\em sector-germ\/} $S_y$ of $S$ if there exists an open set $U$ of $\Delta$ 
such that $U\cap S\subseteq V$ and $y\in U\cap S$.  Similarly if $P$ is a panel of $S$, 
we say that $V$ contains the {\em sector-panel germ} $P_y$ if $y\in U\cap P_y\subset V$.  

We say that the sectors $S_1=w_1(S)$ and $S_2=w_2(S)$ are {\em $i$-adjacent} if $\overline{w_1}=r_1\overline{w_2}$, in $\overline{W}=W/\T$, and we define a {\sl sector-gallery} to be a sequence
$$
S_0,\dots, S_k
$$
of sectors such that $S_{i-1}$ is $j_i$-adjacent to $S_{i}$.  The type of the sector-gallery is $j_1 j_2\dots j_k$.  The {\em sector distance} is then given by $d(S_0,S_k)=r_{j_k}r_{j_{k-1}}\dots r_{j_2} r_{j_1}$.  Given two sector germs $S_y$ and $T_y$ of $\Sig$ (corresponding to sectors $S$ and $T$), the sector germ distance from $S_y$ to $T_y$ is given by $\del(S_y,T_y)=d(S,T)$.  

We note here that if $\Sigma$ were the simplicial affine apartment, a sector germ corresponds to the chamber at the base of the sector.  The distance between two sector germs based at the same point corresponds precisely with the distance between those two chambers in the residue of the point.  On the other hand, the distance between two sectors corresponds to the distance between the parallel classes of those sectors in the apartment at infinity.

We conclude the discussion of the $\Lambda$-apartment with a straightforward lemma that will be used in the sequel
\begin{Lemma} \label{lm:Sec conv}
Let $S$ be a sector of $\Sigma$ based at $x$.  If $S'$ is a subsector of $S$, then the convex hull of $S'$ and $x$ is $S$.
\end{Lemma}

We refer the reader to \cite{B2} for any other definitions needed.

\subsection{Axioms for a $\Lambda$-building}

In this section we will state the various axioms that we will use in this paper.  The first set of six axioms are those in the original definition of an affine $\Lambda$-building.  

Retaining the notation of the previous section, let $\Delta$ be a set, $\Phi$ a spherical root system, $\Sig$ the associated canonical $\Lambda$-apartment, and $\cF$ a set of maps from $\Sig$ to $\Delta$.  An apartment of $\Delta$ is a set $f(\Sigma)$ for some $f\in\cF$. We define the sectors, sector-panels, walls, half-apartments, etc. of $\Delta$ to be the images of the same such under any $f\in\cF$.  

The pair $(\Delta,\cF)$ is an {\em affine $\Lam$-building} if the following conditions are satisfied:
\begin{itemize}
\item[(A1)] Given $f\in \cF$ and $w\in W$, then $f\circ w\in\cF$.

\item[(A2)] (Compatability Axiom) Given $f_1,f_2\in\cF$, if $f_1(\Sig)\cap f_2(\Sig)\ne\emptyset$ then $f_1^{-1}f_2(\Sig)$ is a closed convex set of $\Sig$, and there exists $w\in W$ such that $f_1|_{f_1^{-1}f_2(\Sig)}=f_2\circ w|_{f_1^{-1}f_2}(\Sig)$.

\item[(A3)] Given $x,y\in \Delta$ then there exists $f\in\cF$ such that $x,y\in f(\Sigma)$.  (Note that axioms (A2) and (A3) imply that the metric $d$ on $\Sigma$ extends via $\cF$ to a well defined distance function $d:\Delta\times\Delta\rightarrow\Lambda$.)

\item[(A4)] (Subsector axiom) Given sectors $S_1, S_2\subset \Delta$, then there exists subsectors $S_1'\subset S_1$ and $S_2'\subset S_2$ and $f\in\cF$ such that $S_1'\cup S_2'\subset f(\Sigma)$.  

\item[(A5)] (Retraction axiom) Given $f\in\cF$ and $x\in f(\Sig)$, there exists a retraction $\rho:\Delta\rightarrow f(\Sigma)$ such that $\rho^{-1}(x)=x$.  (Note that this axiom implies the distance function $d$ on $\Delta$ satisfies the triangle inequality (see \cite{B2}).)

\item[(A6)] ($Y$-condition) Given three maps $f_1,f_2,f_3\in\cF$ such that $f_i(\Sigma)\cap f_j(\Sig)$ is a half-apartment for $i\ne j$, then $f_1(\Sig)\cap f_2(\Sig)\cap f_3(\Sig)$ is non-empty.

\end{itemize}

Given this definition of a Affine $\Lambda$-building $(\Del,\cF)$, the set of apartments of $(\Del,\cF)$ is given by
$$
\A=\{f(\Sig)\mid f\in\cF\}.
$$
Abusing notation, we will often refer to the set $\Del$ as satisfying the axioms and assume the set $\A$ of apartments (and maps) is given.

\section{Exchange Axioms}

In this section, we present two new axioms and prove their equivalence to (A6) (given (A1)-(A5)).  Our two new axioms are
\begin{itemize}
\item[(EC)] (Exchange Condition)  Given two maps $f_1,f_2\in\cF$ such that $f_1(\Sig)\cap f_2(\Sig)$ is a half apartment, then there exists a map $f_3\in\cF$ such that $f_3(\Sig)\cap f_j(\Sig)$ is a half apartment for $j=1,2$.  Moreover, $f_3(\Sig)$ is the symmetric difference of $f_1(\Sig)$ and $f_2(\Sig)$ together with the boundary wall of $f_1(\Sig)\cap f_2(\Sig)$.
\end{itemize}
Note that the exchange condition can be restated in ``apartment language'' as: Given two apartments $A_1$ and $A_2$ of $(\Delta,\cF)$ intersecting in a half-apartment $H$ with boundary wall $M$, then $(A_1\oplus A_2)\cup M$ is also an apartment (where $\oplus$ denotes the symmetric difference).
\begin{itemize}
\item[(SE)] (Strong Exchange Condition)  Suppose $f_1\in\cF$ and $S$ is a sector of $(\Delta,\cF)$ such that $P=S\cap f_1(\Sigma)$ is a sector-panel.  Letting $M$ be the wall of $f_1(\Sig)$ containing $P$.  Then there exist $f_2,f_3\in\cF$ such that $f_1(\Sig)\cap f_j(\Sig)$ is a half-apartment and $(M\cup S)\subset f_j(\Sig)$ (for $j=2,3$). 

\end{itemize}

The Strong Exchange Condition can be restated as: Given an apartment $A$ of $\Delta$ and a sector $S$ with a sector-panel of $S$ lying in $A$, then there exists apartments $A_1$ and $A_2$ such that $S\subset (A_1\cap A_2)$ and $A\cap A_j$ is a half-apartment containing the sector-panel of $S$ for $j=1,2$.

Following the notation of \cite[p.571]{B2} and as we did for $\Sigma$, we define two sectors $S_1$ and $S_2$ of $\Delta$ to be parallel if $S_1$ and $S_2$ are at bounded distance from each other.  We similarly define sector panels to be parallel (see \cite{B2} for how we define the type function on sector-panels) if they are at bounded distance. 
We note that axioms (A2) and (A4) imply the sector-distance function on $\Sigma$ extends to a well-defined $\overline{W}$-distance function $d$ on the parallel classes of sectors of $\Delta$.

Let $\Delta$ ($=(\Delta,\cF)$) be an affine $\Lambda$-building of type $W$ and define $\Delta^{\infty}=\{S^{\infty}\mid \mbox{$S$ is a sector of $\Delta$}\}$. We say that $S_1^\infty$ and $S_2^\infty$ are $i$-adjacent if $S_1$ and $S_2$ have parallel sector-panels of type $i$.  For an apartment $A$ of $\Delta$, we define the set $A^\infty=\{S^\infty\mid \mbox{$S$ is a sector of $A$}\}$.  Then if $\Delta$ is an affine $\Lambda$-building of type $W$, the set $A^{\infty}$ is a thin chamber complex of type $\overline{W}$.

\begin{Theorem}{\label{Th:build at inf}} (\cite[Theorem~3.7]{B2})
Let $(\Delta,\cF)$ be a pair satisfying conditions (A1)-(A5).  Then the chamber system $\Delta^\infty$ is a spherical building of type $\overline{W}$ with apartments in one-to-one correspondence with the apartments of $\Delta$.  Similarly, the walls and panels are in one-to-one correspondence with the parallel classes of walls and the parallel classes of sector panels of $\Delta$.
\end{Theorem}

In \cite{B2}, the assumptions are actually that $(\Delta,\cF)$ satisfies condition (A6) also, but an analysis of the proof shows that condition~(A6) is never used.  In fact, the $Y$-condition (A6) is used in \cite{B2} primarily to avoid pathological cases and to force the existence of $\Lambda$-trees associated to the walls and panels at infinity as in \cite{T1}. 

Given Thoerem~\ref{Th:build at inf}, we are now ready to prove the first of our equivalence results.

\begin{prop} \label{Pr:A6 A6'}
Given $(\Del,\cF)$ satisfies conditions~(A1)-(A5), then condition~(A6) is equivalent to condition~(EC).
\end{prop}
\proof
Let us begin by assuming that $(\Del,\cF)$ satisfies (A6), and suppose $A_1=f_1(\Sig)$ and $A_2=f_2(\Sig)$ are two apartments of $\Del$ with $A_1\cap A_2=H$ a half-apartment.  Then $A_1^\infty$ and $A_2^\infty$ are apartments of $\Delta^\infty$ that intersect in a half-apartment.  By spherical building theory, it follows that there exists an apartment $A_3^\infty$ whose chambers are the chambers of $A_1^\infty\oplus A_2^\infty$.  Theorem~\ref{Th:build at inf} now implies that there exists an apartment $A_3$ of $\Delta$ corresponding to $A_3^{\infty}$.  Consequently by (A6) $A_1\cap A_2\cap A_3$ is non-empty.  Since $A_1^{\infty}\cap A_3^{\infty}$ is a half-apartment and $A_1\cap A_3$ is convex by condition~(A2), it follows that $A_1\cap A_3$ is a half-apartment.  Similarly $A_2\cap A_3$ is a half-apartment.  Condition~(A6) now implies that $A_1\cap A_2\cap A_3$ contains some element $x\in\Delta$.  Since $x\in H$ and $A_3^\infty$ contains the chambers of $A_1^\infty$ not in $A_2^\infty$, it follows that $A_3$ contains $A_1-{H}$.  Similarly $A_2 - {H}\subset A_3$.  By convexity $\bound{H}\subseteq A_3$.  But now the convexity of $A_3$ implies that $x\in \bound H$ (the boundary of $H$) as otherwise the wall parallel to $\bound H$ through $x$ would not separate points of $A_1\cap A_3$ and $A_2\cap A_3$.  This implies that (EC) holds.

Now assume (A1)-(A5) and (EC) are all satisfied, and let $A_1$, $A_2$, and $A_3$ be half apartments of $\Delta$ such that any two intersect in a half-apartment.  By way of contradiction, suppose $A_1\cap A_2\cap A_3=\emptyset$.  Let $H_{ij}=A_i\cap A_j$ for $i,j\in\{1,2,3\}$.  Since $H_{1,2}\cap H_{1,3}=\emptyset$, it follows that if $H$ is a half-apartment of $A_1$ with $H_{1,2}\cap H\subseteq\bound H$, then $H_{1,3}\cap H$ is again a half-apartment.  Now (EC) implies that there exists an apartment $A_4$ such that $A_4=(A_1\oplus A_2) \cup \bound H_{1,2}$.  Note that $H_{1,3}\subseteq A_4$, so that $A_4^{\infty}$ consists of the same sectors as $A_3^{\infty}$.  However, by Theorem~\ref{Th:build at inf}, the apartments of $\Delta$ are in one-to-one correspondence with the apartments of $\Delta^{\infty}$.  Therefore, $A_3=A_4$. 
\qed

We now prove the similar result for condition (SE).
\begin{prop} \label{Pr:A6' A6''}
Given $(\Del,\cF)$ satisfies conditions~(A1)-(A5), then condition~(EC) is equivalent to condition~(SE).
\end{prop}
\proof  Suppose $(\Del,\cF)$ satisfies conditions~(A1)-(A5) and (EC).  Suppose $A_1$ is an apartment and $S$ is a sector such that $S\cap A_1$ is a sector panel of $S$.  Again, by Theorem~\ref{Th:build at inf}, in $\Del^{\infty}$, $A_1^\infty\cap S^\infty$ is a panel.  Therefore, for the building at infinity, there is an apartment $A_2^\infty$ of $\Delta^\infty$ such that $S^\infty\in A_2^\infty$, and $A_1^\infty\cap A_2^{\infty}$ is a half apartment.  Let $A_2$ be the corresponding apartment of $\Delta$.  Since $A_1^\infty\cap A_2^{\infty}$ is a half-apartment of $\Delta^\infty$, it follows that $A_1\cap A_2$ is a half-apartment of $\Delta$.  We now apply condition~(EC) to conclude the argument.

Conversely, suppose $(\Del,\cF)$ satisfies conditions~(A1)-(A5) and~(SE), and let $A_1$ and $A_2$ be apartments of $\Delta$ intersecting in a half-apartment $H$.  Let $S$ be a sector of $A_2$ such that $S\cap A_1$ is a sector-panel $P$ of $S$, and let $M$ be the wall of $A_1$ containing $P$.  By (EC), there exists an apartment $A_3$ containing $M$ such that $A_1\cap A_3$ is a half-apartment and $A_2\cap A_3$ is a half-apartment (as $A_2$ must be one of the apartments guaranteed by (EC)) containing $M\cup S$. By convexity, it follows that $A_3=(A_1\oplus A_2)\cup M$ as desired.
\qed

We have now shown that the $Y$-condition can be replaced by either of the exchange axioms.   

\section{Main Results}

In this section we prove our two main theorems.  The first theorem provides a general proof for affine buildings ($\Lambda$- or otherwise) that for any two sector-germs of $\Delta$ are contained in a common apartment.  That is, our result generalizes the work of Anne Parreau \cite{P} on Bruhat-Tits buildings, although our proof takes a different approach.  

The second theorem states $\Delta$ is an affine $\Lambda$-building if $\Delta$ satisfies axioms (A1), (A2), (A3), (A4), and (SE).  The proof of the second theorem uses a slightly stronger result than the first theorem.  

To begin with we need a preparatory lemma.
\begin{Lemma} \label{Lm:secgm sec}
Let $\Delta$ satisfy conditions (A1)-(A4) and (SE).  If $S$ and $T$ are sectors of $\Delta$  based at $y$, there exists an apartment $A$ of $\Delta$ such that $S_y \cup T\subset A$.  Moreover, $\ell(\delta(S_y,T_y))\le\ell(d(S^{\infty},T^{\infty}))$ and equality holds if and only if $S$ and $T$ are contained in a common apartment.
\end{Lemma}
By axiom (A4) there exists an apartment $A'$ containing subsectors $S'$ of $S$ and $T'$ of $T$, and $d(S',T')=d(S,T)$.  Consider the sector-gallery
$$
S'=S_0,\dots,S_n=T',
$$
and let $A_0$ be an apartment containing $S$ (and hence $S'$).  Let $j$ be minimal such that $S_{j+1}$ contains no subsector in $A_0$.  If $j=n$ exists (i.e., $T'$ has a subsector $T''$ contained in $A_0$), then $T''$ is a sector of $A_0$, and as $y\in A_0$, by convexity it follows that $T\subset A_0$ and there is nothing to prove.  We will induct on $n-j$.  The basis step having been proven, assume $S_{j+1}$ has no subsector contained in $A_0$ but $S_0,\dots,S_j$ all have subsectors in $A_0$.  In this case, there exists a sector $S_{j+1}'$ parallel to $S_{j+1}$ (in $A'$) such that $S_j'\cap A_0$ is a sector-panel (parallel to a sector panel of $S_j$).  By condition (SE) there exists an apartment $A_{j+1}$ containing $S_{j+1}'$ and the sector germ $S_y$ (since for any wall $S_y$ must lie on one side or the other of the wall).  If $S$ is contained in $A_{j+1}$, then we replace $A_0$ with $A_{j+1}$ and by induction on $n-j$ we have the result.  On the other hand, if $S\not\subset A_{j+1}$, let $S''$ be the sector of $A_{j+1}$ with sector germ $S''_y=S_y$.  Then $$
\ell(d(T',S_{j+1}))=\ell(d(T,S_{j+1}))-1.
$$
Moreover, considering the case of $S''$ and $T$, together with $A_{j+1}$ as our new $A_0$, by induction there exists an apartment $A$ containing $S''_y$ and $T$, with $\ell(\delta(S''_y,T_y))\le \ell(d(S''^\infty,T^\infty)$.  However, $S_y=S''_y$ and 
$$
\ell(d(S''^\infty,T^\infty))\le\ell(d(S^\infty,T^\infty))-1.
$$  
Hence
$$
\ell(\del(S_y,T_y))\le\ell(d(S,T))
$$
as desired.  Note that if equality holds, then in each case, the apartment $A_j$ contains $S$ (where we take $A_j$ as the apartment containing $S_j$ and $S_y$ in the proof), and in particular, $A_n$ contains both $S$ and $T$ as desired.
\qed

\begin{cor} \label{Lm:opp sec}
If $S$ and $T$ are sectors of $\Delta$ based at $y$, and $\delta(S_y,T_y)$ is maximal, then $S$ and $T$ are contained in a common apartment.
\end{cor}
\proof  Since $\delta(S_y,T_y)$ is maximal, Lemma~{\ref{Lm:secgm sec}} implies that $\ell(\delta(S_y,T_y))=\ell(d(S,T))$.  However, in this case the lemma implies the existence of an apartment $A$ containing $S$ and $T$. 
\qed

From here we can now prove the first of our two theorems.
\begin{Theorem} \label{Th:sec sec}
Suppose $\Delta$ satisfies (A1), (A2), (A3), (A4), and (SE).  Let $S$ and $T$ be sectors of $\Delta$ based at $x$ and $y$ respectively.  Then there exists an apartment $A$ of $\Delta$ containing $S_x$ and $T_y$. That is, two sector germs are contained in a common apartment.
\end{Theorem}
\proof
We begin by showing that $S_x$ and $y$ are contained in a common apartment $B$.  By (A2), there exists an apartment $A'$ containing $x$ and $y$.  Let $T'$ be a sector of $A'$ based at $x$ containing $y$.  By Lemma~\ref{Lm:secgm sec} there exists an apartment $B$ of $\Delta$ containing $S_x$ and $T'$.  Hence $B$ contains $S_x$ and $y$.  Take $S'$ to be the sector of $B$ based at $y$ containing $S_x$.  Again by Lemma~\ref{Lm:secgm sec}, there exists an apartment $A$ containing $T_y$ and $S'$.  Since $S_x\subset S'$, it follows that $A$ contains $T_y$ and $S_x$ as desired.
\qed

We now wish to show that condition (SE) can replace conditions (A5) and (A6) in the definition of an affine $\Lambda$-building.  Since we have already shown that (SE) is satisfied by an affine $\Lambda$-building $\Delta$, and that (A6) can be replaced by (SE), it remains to show that (A1)--(A4) and (SE) together imply the sector-retraction condition (A5).  Our argument will proceed along the following lines.  First we use Lemma~\ref{Lm:secgm sec} to define the retraction $\rho_{S_x,A}$ for an apartment $A$ and a sector-germ $S_x\subset A$ and show that this retraction preserves distances on sets $X$ contained in a common apartment with $S_x$.  We then will show that (SE) implies a slightly stronger exchange condition, namely that if $S_x$ is a sector germ having a sector-panel germ in an apartment $A$, then there exist apartments $A'$ and $A''$ such that $S_x\subset A'\cap A''$ and $A\subseteq A'\cup A''$.  Given this condition, we will then be able to show that given an apartment $B$, a point $y\in B$, and a sector-germ $S_x$, there exists a sector $T$ containing $y$ based at $x$ such that $B$ contains a subsector of $T$.  Since any two sectors based at $x$ that contain a common subsector must be equal, it follows that $A$ is contained in a finite union of closed convex subsets $X_1$, \dots, $X_n$, each of which are contained in an apartment with $S_x$.  Consequently, the retraction $\rho_{S_x,A}$ preserves distances on each of the $X_i$s.  Given any pair of points $y$ and $z$ of $B$, it follows that we can find points 
$$
y=y_0,y_1,\dots,y_t=z
$$
such that $y_{i-1},y_i\in X_{j_i}$ and $y_i$ is in the convex hull of $y_{i-1}$ and $y_{i+1}$.   As a result, $\rho_{S_x,A}$ will preserve the distance between consecutive $y_i$s, and using the triangle inequality on $A$, it will follow that 
\begin{eqnarray*}
d(\rho_{S_x,A}(y),\rho_{S_x,A}(z)) & \le & \sum_{i=1}^n d(\rho_{S_x,A}(y_{i-1}),\rho_{S_x,A}(y_i)) \\
& = & \sum_{i=1}^n d(y_{i-1},y_i) =d(y,z)
\end{eqnarray*}

Let $A$ be an apartment of $\Delta$, $x\in A$, and $S$ a sector of $A$ based at $x$.  Let $S_x$ be the sector germ of $S$.  We define $\rho_{S_x,A}$ as follows:  For $y\in \Delta$, Theorem~\ref{Th:sec sec} there exists an apartment $B$ containing $y$ and $S_x$.  Let $f\in\F$ be such that $f(\Sig)=B$, $g\in\F$ be such that $g(\Sig)=A$.  Since $S_x\in A\cap B$, by~(A2) there exists $w\in W$ such that $g(w(f^{-1}(S_x)))=S_x$.  Moreover, as the identity is the only element of $W$ fixing a sector-germ, $w$ is unique.  Define $\rho_{S_x,A}(y)=g(w(f^{-1}(y)))$.  By the compatibility condition, if $B'$ is another apartment of $\Delta$ containing $S_x$ and $y$ with $f'(\Sig)=B'$, it follows that there exists $w'\in W$ with $g(w'(f'^{-1}(S_x)))=S_x$.  This implies that 
$$
w'(f'^{-1}(S_x))=w(f^{-1}(S_x)).
$$
There also exists an element $w''\in W$ such that $w''(f'^{-1}(S_x))=f^{-1}(S_x)$.  Sine the stabilizer of a sector germ is trivial, $w''=w^{-1}w'$.  But $w''(f'^{-1}(y))=f^{-1}(y)$ by the compatibility condition (A2).  Therefore, $w'f'^{-1}(y)=wf^{-1}(y)$, implying $g(w'f'^{-1}(y))=g(wf^{-1}(y))$, and thus $\rho_{S_x,A}$ is well-defined.  Since the map $w\in W$ preserves the distance on $\Sig$, it follows that $d(y,z)=d(\rho_{S_x,A}(y),\rho_{S_x,A}(z))$ for all $y,z\in\Delta$ such that $y$, $z$, and $S_x$ are contained in a common apartment.  

We next show an even stronger exchange condition is satisfied.
\begin{Lemma}\label{Lm:sup exch}
Suppose $\Delta$ satisfies conditions (A1), (A2), (A3), (A4), and (SE).  Let $A$ be an apartment of $\Delta$ and $S_x$ be a sector germ of $\Delta$ such that $S_x\cap A$ is a sector-panel germ $P_x$. Then there exists apartments $A'$ and $A''$ such that $S_x\in A'\cap A''$ and $A\subset A'\cup A''$.
\end{Lemma}
\proof
By the strong exchange condition (SE) it suffices to show that there is a sector $S'$ of $\Delta$ intersecting $A$ in a sector-panel such that $S'_x=S_x$.  Let $T$ be a sector of $A$ based at $x$ having a sector-panel containing the sector panel germ of $S_x$.  (That is, if $P'$ is the sector-panel of $A$ having sector panel germ $P_x$, take $T$ to be a sector having a sector-panel $P'$.)  By Lemma~\ref{Lm:secgm sec} there exists an apartment $B$ containing $T$ and $S_x$.  Let $S'$ be the sector of $B$ having sector-panel germ $S_x$.  Then $S'$ has sector-panel $P'$.  Moreover, by convexity, if $S'\cap A\ne P'$, then $S_x=S'_x\subset A$ contrary to our hypothesis. Therefore $S'\cap A=P'$ and by (SE) there exists apartments $A'$ and $A''$ such that $S'\subset A'\cap A''$ and $A\subset A'\cup A''$.
\qed

This exchange condition allows us to work with sector-germs based at a common point, much as in the simplicial buildings case one works with chambers in a spherical residue. 

\begin{prop} \label{Pr:opps}
Suppose $\Delta$ satisfies conditions (A1), (A2), (A3), (A4), and (SE).  Let $S_x$ be a sector germ contained in an apartment $A$ of $\Delta$, and $B$ be another apartment of $\Delta$.  Then for every point $y\in B$ there exists a sector $T$ of $B$ such that 
\begin{enumerate}
\item There exists a sector $T'$ based at $x$ parallel to $T$ containing $y$, and
\item There exists an apartment $A'$ of $\Delta$ containing $T$ and $S_x$.
\end{enumerate}
\end{prop}
\proof  By Theorem~\ref{Th:sec sec}, for every sector $T$ of $B$ based at $y$, there exists an apartment $B'$ of $\Delta$ containing $T_y$ and $S_x$ let $S'$ denote the sector of $B'$ based at $y$ containing $S_x$.  For $y\in B$, choose $T$ such that $\ell(\delta(T_y,S'_y))$ is maximal.  If $T_y$ and $S'_y$ are not opposite (that is $\delta(T_y,S'_y)$ is not the longest element of $\overline{W}$) then let $P_y$ be a sector-panel germ of $T_y$ such that the wall $M$ of $B'$ through $P_y$ does not separate $T_y$ and $S'$.  In the apartment $B$ there exists a sector $R$ such that $R_y$ shares $P_y$ with $T_y$.  Moreover, since $T_y$ and $S'$ lie on the same side of $M$, by convexity the apartment $B''$ containing $S'$ and $R_y$ guaranteed by Lemma~\ref{Lm:sup exch} also contains $T_y$.  In $B''$, we then have $\ell(\delta(R_y,S'_y))=\ell(\delta(T_y,S'_y)+1$, contradicting the choice of $T$.  Hence we may assume that $T_y$ and $S'_y$ are opposite.  By Corollary~\ref{Lm:opp sec}, there exists an apartment $A'$ of $\Delta$ containing $S'$ and $T$.  But $S_x\subset S'$, so that $S_x\subset A'$.  Moreover, since $A'$ contains $T$, take $T'$ be the sector based at $x$ parallel to $T$ (in $A'$).  Since $T$ and $S'$ were opposite sectors and $x\in T$, it follows that $y\in T'$, completing the proof of the proposition.
\qed

Note that by convexity, there is at most one sector of $\Delta$ parallel to $T$ based at a point $x$.  Consequently, Proposition~\ref{Pr:opps} implies that there exist finitely many sectors $S_1$, $S_2$, \dots, $S_n$ of $\Delta$ each based at $x$, such that $S_x$ is contained in a common apartment $A_i$ with each of the $S_i$.  

\begin{Corollary}\label{Lm:fin cover}
Suppose $\Del$ satisfies conditions (A1), (A2), (A3), (A4), and (SE).  Let $S_x$ be a sector germ contained in an apartment $A$ of $\Delta$ and $B$ be an apartment of $\Del$.  Then there exists closed convex sets $X_1$, \dots, $X_n$ of $B$ such that 
\begin{enumerate} 
\item $B=X_1\cup\dots\cup X_n$ and
\item Each $X_i$ lies in a common apartment with $S_x$.
\end{enumerate}
\end{Corollary}
\proof  Let $S_1$, \dots, $S_n$ be the sectors from the above paragraph, and set $X_i=S_i\cap B$.
\qed

We are now ready to prove the map $\rho_{S_x,A}$ is a retraction.  
\begin{prop}
Suppose $\Del$ satisfies conditions (A1), (A2), (A3), (A4), and (SE).  Let $S_x$ be a sector germ contained in an apartment $A$ of $\Delta$. Then the map $\rho_{S_x,A}$ is a retraction of $\Del$ onto $A$.
\end{prop}
Let $\rho=\rho_{S_x,A}$.  We have already shown that $\rho:\Del\rightarrow A$ is well-defined, and $\rho|_A$ is the identity by definition.  It remains to see that $\rho$ diminishes distance.  That is, to show for all $y,z\in\Del$ that
$$
d(\rho(y),\rho(z))\le d(y,z).
$$
By definition, if there is an apartment $B$ of $\Del$ containing $S_x$, $y$, and $z$, then $d(y,z)=d(\rho(y),\rho(z))$, so the result holds true.  Now, suppose $y$ and $z$ are arbitrary.  By (A2) there exists an apartment $B$ containing $y$ and $z$.  By Corollary~\ref{Lm:fin cover} there exists closed convex sets $X_1$, \dots, $X_n$ such that $B=\cup_{i=1}^n X_i$ and each $X_i$ is contained in a common apartment with $S_x$.  Since each $X_i$ is convex and closed, there exists a sequence of points
$$
y=y_0,y_1,\dots,y_k=z
$$
such that $y_{i-1},y_i\in X_{j_i}$ for some $j_1$, \dots, $j_k$ and $y_i$ is in the convex hull of $y_{i-1}$ and $y_{i+1}$ for $i=1,\dots,k-1$.  Then
\begin{eqnarray*}
d(y,z)&=&\sum_{i=1}^k d(y_{i-1},y_i) \\
  & = & \sum_{i=1}^k d(\rho(y_{i-1}),\rho(y_i)) \\
  & \ge & d(\rho(y_0),\rho(y_k)) \\
  & = & d(\rho(y),\rho(z),
\end{eqnarray*}
where we use the triangle inequality for $d$ restricted to $\Sigma$ in the next to the last step.  Thus, $\rho$ is distance diminishing and hence a retraction with the required properties.
\qed

Summarizing, we have proven
\begin{Theorem} \label{Th:main}
Suppose $(\Del,\cF)$ satisfies conditions (A1), (A2), (A3), and (A4).  Then conditions (A5) and (A6) together are equivalent to condition (SE).  In particular, if $(\Del,\cF)$ satisfies conditions (A1), (A2), (A3), (A4), and (SE) then $(\Del,\cF)$ is an affine $\Lam$-building.
\end{Theorem}
\qed.

\ignore{\begin{thebibliography}[BT]

\bibitem[B1]{B1} Bennett, Curtis, Affine $\Lambda$-buildings, Ph.D. Thesis, University of 
Chicago, 1990.

\bibitem[B2]{B2} Bennett, Curtis, Affine $\Lambda$-buildings I, {\it Proceedings of the 
London Math Society (3) {\bf 68} (1994), no. 3, 244--267.

\bibitem[BT]{BT} Bruhat, F., Tits, J., Groupes r\'eductifs sur un corps local, {\it Inst. Hautes 
\'Etudes Sci. Publ. Math.}, no. 41 (1972), 5--251.

\bibitem[DS]{DS} Delgado, A., Stellmacher, B.  Weak $(B,N)$-pairs of rank~2, {\it Groups 
and Graphs: New Results and Methods}, DMV Seminar, 6, Birkh\"auser, Basel, 1985, 
59--244.

\bibitem[KT]{KT} Kramer, L., Tent, K., Asymptotic cones and ultrapowers of Lie groups, {\it 
Bulletin of Symbolic Logic} {\bf 10} (2004), no. 2, 175-185.

\bibitem[MS]{MS} Morgan, J., Shalen, P., Valuations, trees, and degenerations of hyperbolic 
structures, I., {\it Ann. Of Math. (2)} {\BF 120} (1984), no. 3, 401-476.

\bibitem[P]{P} Parreau, A., Immeubles affines: construction par les normes et \'etude des isom\'etries, {\it Crystallographic groups and their generalizations (Kortrijk 1999)}, 266--302, {\it Contemporary Mathematics}, {\bf 262}, American Mathematical Society, Providence, RI, 2000.

\bibitem[T]{T1} Tits, J., Immeubles de type affine, {\it Buildings and the geometry of 
diagrams (Como, 1984)}, 159-190, Lecture Notes in Mathematics, 1181, Springer-
Verlag, Berline, 1986.}}

\end{document}